\documentclass[12pt]{article}
\usepackage{latexsym}
\usepackage{amsmath}
\usepackage{amssymb}
\usepackage{amsfonts}
\usepackage{chngcntr}

\numberwithin{equation}{section}

\newtheorem{theorem}{Theorem}[section]

\newtheorem{definition}{Definition}[section]
\newcommand{\eproof}{{\mbox{\ }~\hfill
\mbox{\large $\Box$} \par \vskip 10pt}}

\newcommand{\R}{{\mathbb R}}

\renewcommand{\div}{{\rm div}}

\title{Doubling inequalities for the Lam\'e system with rough coefficients}
\author{Herbert Koch\thanks{Mathematisches Institut, Universit\"at Bonn, Endenicher Allee
60, D-53115 Bonn, Germany. Partially supported by the DFG through SFB 1060. Email:koch@math.uni-bonn.de}\qquad
Ching-Lung Lin\thanks{Department of Mathematics and Research Center for Theoretical Sciences, NCTS, National Cheng-Kung University, Tainan 701, Taiwan.  Partially supported by the Ministry of Science and Technology of Taiwan. Email:
cllin2@mail.ncku.edu.tw}
\qquad Jenn-Nan Wang\thanks{Institute of Applied Mathematical Sciences, NCTS, National Taiwan
University, Taipei 106, Taiwan. Partially supported by MOST102-2115-M-002-009-MY3. Email:jnwang@math.ntu.edu.tw}}

\date{}

\begin{document}
\maketitle

\begin{abstract}
In this paper we study the local behavior of a solution to the
Lam\'e system when the Lam\'e coefficients $\lambda$ and $\mu$ satisfy that $\mu$ is Lipschitz and
$\lambda$ is essentially bounded in dimension $n\ge 2$. One of the main results is the \emph{local} doubling inequality
for the solution of the Lam\'e system. This is a quantitative
estimate of the strong unique continuation property. Our proof
relies on Carleman estimates with carefully chosen weights.
Furthermore, we also prove the \emph{global} doubling inequality,
which is useful in some inverse problems.
\end{abstract}

\section{Introduction}\label{sec1}
\setcounter{equation}{0}

Let $\Omega$ be an open connected  subset of $\R^n$ with $n\ge 2$.
Without loss of generality, we assume $0\in\Omega$. Let $\mu(x)\in C^{0,1}(\Omega)$ and
$\lambda(x),\rho(x)\in L^{\infty}(\Omega)$ satisfy
\begin{equation}\label{1.1}
\begin{cases}
\mu(x)\geq\delta_0,\quad\quad
\lambda(x)+2\mu(x)\geq\delta_0\quad\forall\
x\in\Omega,\\
\|\mu\|_{C^{0,1}(\Omega)}+\|\lambda\|_{L^{\infty}(\Omega)}\leq
M_0,\quad \|\rho\|_{L^{\infty}(\Omega)}\le M_0
\end{cases}
\end{equation}
with positive constants $\delta_0, M_0$, where we define
$$\|f\|_{C^{0,1}(\Omega)}=\|f\|_{L^{\infty}(\Omega)}+\|\nabla
f\|_{L^{\infty}(\Omega)}.$$
The isotropic elasticity system is given by
\begin{equation}\label{1.2}
\text{div}(\mu(\nabla u+(\nabla
u)^t))+\nabla(\lambda\text{div}u)+\rho u=0\quad\text{in}\ \Omega,
\end{equation}
where $u=(u_1,u_2,\cdots,u_n)^t$ is the displacement vector and $(\nabla
u)_{jk}=\partial_ku_j$ for $j,k=1,2,\cdots,n$. If $\rho=0$, \eqref{1.2} represents the displacement
equation of equilibrium.

Under the assumptions \eqref{1.1}, the qualitative strong unique
continuation property for \eqref{1.2} were recently proved by
Nakamura, Uhlmann and the second and third authors \cite{lin6}, i.e., if $u\in H^1(\Omega)$
solves \eqref{1.2} and satisfies that for any $N\in{\mathbb N}$,
there exists a constant $C_N$ such that
$$\int_{B_r}|u|^2\le C_Nr^N\quad\forall\ r\ \text{sufficiently small},$$ then $u\equiv 0$
in $\Omega$. In fact, in \cite{lin6}, we derived a quantitative
estimate on the vanishing order of any nontrivial solution to
\eqref{1.2}. The derivation relies on the \emph{optimal} three-ball
inequalities (see \cite{lin6} for details).

Another quantitative estimate of the strong unique continuation
property is the doubling inequality. When $\lambda,\mu\in C^{1,1}$
and $\rho=0$, doubling inequalities for \eqref{1.2} in the form
\begin{equation*}
\int_{B_{2r}}|u|^2+|\div u|^2\le K\int_{B_r}|u|^2+|\div u|^2.
\end{equation*}
were derived in \cite{almo} based on the frequency function method
developed in \cite{gl1} and \cite{gl2}. To apply quantitative
estimates of the strong unique continuation property to certain
inverse problems for the elasticity, it is desirable to derive a
doubling inequality containing $|u|^2$ only \cite{amr}, i.e.,
\begin{equation}\label{2doub}
\int_{B_{2r}}|u|^2\le K\int_{B_r}|u|^2.
\end{equation}
Indeed, \eqref{2doub} for the Lam\'e system with $C^{1,1}$
coefficients was proved in \cite{amr}. However, as mentioned in
\cite{amrv}, the proof given there contains a gap.

In \cite{amrv}, the authors proved doubling inequalities of the form
\eqref{2doub} when $\lambda,\mu\in C^{2,1}$ (also $\rho=0$).
Moreover, these inequalities depend on global properties of the
solution. A key observation in \cite{amrv} is that for
$\lambda,\mu\in C^{2,1}$, the Lam\'e system can be transformed into
a fourth order system for $u$ having $\Delta^2$ as the leading part
and essentially bounded coefficients in the lower orders. For this
fourth order system, three-sphere inequalities and local doubling
inequalities were derived in \cite{lin3}. Using these inequalities,
global doubling inequalities for \eqref{1.2} were then obtained.

The aim of this paper is to establish doubling inequalities of the
form \eqref{2doub} for \eqref{1.2} when $\mu\in C^{0,1}$ and
$\lambda, \rho\in L^{\infty}$. Our result provides a positive answer to
the open problem posed in \cite{amrv} about the doubling inequality
for \eqref{1.2} with less regular coefficients. The ideas of our
proof originate from our series papers on proving quantitative
uniqueness for elliptic equations or systems by the method of
Carleman estimates \cite{lin4}, \cite{lin5}, and \cite{lin6}. In particular, we will use the reduced
system derived in \cite{lin6} (see Section~2 below).

We now state main results of the paper. Their proofs will be given
in the subsequent sections. Assume that there exists
$0<\tilde{R}_0\le 1$ such that $B_{\tilde{R}_0}\subset\Omega$.
Hereafter $B_r$ denotes an open ball of radius $r>0$ centered at the
origin.

\begin{theorem}\label{thm1.1}
There exists $C>0$ depending on $n$, $M_0$ and $\delta_0$
so that the following is true.
If  $R>0$ with $3R\le\tilde R_0$, $u\in H^1_{loc}(B_R)$ is a nonzero solution to \eqref{1.2} and
\[  m = -  \ln \Big( \frac{ \Vert u\Vert_{L^2(B_{R}\backslash B_{R/2})}}{\Vert u\Vert_{L^2(B_{2R}\backslash B_{R})}}\Big)\]
then
\begin{equation}\label{1.8}
\Vert u\Vert_{L^2(B_r)}\ge C(r/R)^{Cm}   \Vert u \Vert_{L^2(B_{2R}\backslash B_R)}\quad\text{for all} \ r\le R.
\end{equation}
\end{theorem}

\begin{theorem}\label{thm1.2}
There exists a positive constant $\tilde C$ depending only on
$n$, $M_0$ and $\delta_0$ such that the following is true.
If $u\in H^1_{loc}(B_R)$ is a nonzero solution to \eqref{1.2}
 then
\begin{equation}\label{1.9}
\Vert u\Vert_{L^2(B_{2r}(x_0))}\leq \tilde{C} e^{\tilde C m} \Vert u\Vert_{L^2(B_r(x_0))},
\end{equation}
whenever $B_{2r}(x_0) \subset B_{R/2}$. Here $m$ is the constant of
Theorem~\ref{thm1.1}.
\end{theorem}
Theorem~\ref{thm1.1} and \ref{thm1.2} will be proved together. The
estimate \eqref{1.9} is called local doubling inequalities. Global
doubling inequalities in which constants depend on the global
property of solution will be proved in Section~\ref{sec4}.

\section{Reduced system}\label{sect2}
\setcounter{equation}{0}

We now recall the reduced system derived in \eqref{1.2}. This is a
crucial step in our approach. Let us write \eqref{1.2} into a
non-divergence form:
\begin{equation}\label{1.4}
\mu\Delta u+\nabla((\lambda+\mu)\ \div u)+(\nabla u+(\nabla
u)^t)\nabla\mu-\div u\nabla\mu+\rho u =0.
\end{equation}
Dividing \eqref{1.4} by $\mu$ yields
\begin{eqnarray}\label{1.5}
&&\Delta u+\frac{1}{\mu}\nabla((\lambda+\mu)\ \div u)+(\nabla u+(\nabla u)^t)\frac{\nabla\mu}{\mu}-\div u\frac{\nabla\mu}{\mu}+\frac{\rho}{\mu} u\notag\\
&=&\Delta u+\nabla(a(x)p)+G\notag\\
&=&0,
\end{eqnarray}
where $$a(x)=\frac{\lambda+\mu}{\lambda+2\mu}\in
L^{\infty}(\Omega),\quad p=\frac{\lambda+2\mu}{\mu}\ \div u$$ and
$$G=(\nabla u+(\nabla u)^t)\frac{\nabla\mu}{\mu}-\div u(\frac{\nabla\mu}{\mu}+(\lambda+\mu)\nabla(\frac{1}{\mu}))+\frac{\rho}{\mu} u.$$
Taking the divergence on \eqref{1.5} gives
\begin{equation}\label{1.6}
\Delta p+\div G=0.
\end{equation}
Our reduced system now consists of \eqref{1.5} and \eqref{1.6}. It
follows easily from \eqref{1.6} that if $u\in H^1_{loc}(\Omega)$,
then $p\in H^1_{loc}(\Omega)$.

\begin{equation}\label{up}
\begin{cases}
\Delta u+\nabla(a(x)p)+G(x,u)=0,\\
\Delta p+\div G(x,u)=0.
\end{cases}
\end{equation}
Note that system \eqref{up} is not decoupled. We will use \eqref{up}
to prove our theorems.

\section{Proofs of Theorem \ref{thm1.1} and \ref{thm1.2}}\label{sec3}
\setcounter{equation}{0}

This section is devoted to the proofs of Theorem~\ref{thm1.1} and \ref{thm1.2}. The proofs rely on a suitable Carleman estimate proved in \cite{kt01}. To state the estimate, we consider the equation
\begin{equation}\label{eq3-1}
\Delta u+\nabla f=g\quad\mbox{in}\quad\R^n.
\end{equation}
We consider $t>0$. Given $\tau\gg 1$, let $h(t)$ be a convex function satisfying
\begin{equation}\label{weight}
\left\{
\begin{aligned}
&h'\sim \tau,\ \mbox{i.e.},\ \exists\, C>1,\ C^{-1}\tau\le h'\le C\tau,\\
&\mbox{dist}(2h',{\mathbb Z})+h''\gtrsim 1.
\end{aligned}
\right.
\end{equation}
Here and in the sequel, the notation $X\lesssim Y$ or $X\gtrsim Y$
means that $X\le CY$ or $X\ge CY$ with some constant $C$ depending
only on $n$, $M_0$ and $\delta_0$. We further assume that $h$
satisfies that for any $C>0$ there exists $R_0>0$ such that
\begin{equation}\label{eq333}
C|x|\tau\le (1+h''(-\ln |x|))
\end{equation}
for all $\tau$ and $|x|\le R_0$. Given $R>0$ $h(-\ln(R_0x/R)) $
satisfies \ref{eq333} for $|x|\le R$.

 For our purpose, in addition to \eqref{weight}, we also require $h-t-\frac 12\ln(1+h'')$ to
satisfy \eqref{weight}. The existence of such weight function $h$ can be found in \cite[Section~6]{kt01}. We will give a more explicit construction of $h$ in appendix.
\begin{theorem}\label{lem2.3}
Assume that a convex $h$ satisfies \eqref{weight} and is evaluated at $-\ln|x|$. For smooth functions $u$, $f$, $g$ satisfying \eqref{eq3-1} and are supported in $B_1(0)\setminus\{0\}$, we have that
\begin{equation} \label{carleman1}
\begin{split}
&\tau \Vert   |x|^{-2} (1+h'')^\frac12 e^h  u \Vert + \Vert  |x|^{-1} (1+h'')^\frac12 e^h  \nabla u \Vert\\
&\lesssim \tau  \Vert |x|^{-1}  e^h f \Vert +\Vert e^h g\Vert ,
\end{split}
\end{equation}
where $\|\cdot\|=\|\cdot\|_{L^2(\R^n)}$.
\end{theorem}
Theorem~\ref{lem2.3} can be proved by adopting arguments of Proposition~4.1 and 4.2 in \cite{kt01} (see also \cite[Proposition~5.1]{ck10}) .  It can be also proved by modifying the method in \cite{lin4}. Here we give a sketch of proof. 

\medskip\noindent{\it Proof}. We first observe that the estimate is equivalent to the some estimates 
for functions on $\R^n \backslash \{0\}$ under appropriate decay conditions 
of the solutions at $0$ and $\infty$. This is seen by truncating, and taking an obvious limit. 

We begin with an elliptic  reduction and consider the equation 
\[   -\Delta w +  K \tau^2|x|^{-2}  w = \nabla f   \]
with a fixed large positive constant $K$. The quadratic form 
\[  \int |\nabla w|^2 dx + K\tau^2 \int |w|^2 |x|^{-2} dx \] 
is an inner product  and the Riesz representation theorem  ensures that there is a unique 
solution. We claim that 
\begin{equation} \label{weighted} 
 \int e^{2h(-\ln|x|)} (|\nabla w|^2 +\tau^2 |w|^2/|x|^2) dx
\le  C\int e^{2h(-\ln|x|)}| f |^2dx
\end{equation} 
for all $h$ satisfying the first condition of \eqref{weight}. It suffices to consider bounded functions 
$h$ and the inequality follows by multiplying by $e^{2h(-\ln|x|)}w$ and integrating by parts. 
Moreover $w$ decays fast as $x\to  0 $ or $|x|\to \infty$ which we see by choosing $h$ growing fast and linearly at $\pm \infty$. 

We make the ansatz 
\[ u = v+w \] 
where 
\[ \Delta v = -K\tau^2|x|^{-2} w  +g \] 
and the full estimate \eqref{carleman1} follows once we prove the estimate for $f=0$ and apply    it to $v$. Without loss of generality we assume $f=0$ in the sequel and prove \eqref{carleman1}.

To prove it in this case , we introduce polar
coordinates in ${\mathbb R}^n \backslash {\{0\}}$ by setting $x=r
\omega$, with $r=|x|$, $\omega=(\omega_1,\cdots,\omega_n)\in S^{n-1}$. 
Using the   new coordinate $t=-\log r$, we obtain that
\[  |x|^{\frac{2+n}2} \Delta (|x|^{\frac{2-n}2} u) 
= u_{tt} + \Delta_{S^{n-1}} u - \left( \frac{n-2}2\right)^2 u = 
 e^{-\frac{n+2}{2}t}   g(e^{-t}  \omega) . \]
We can diagonalize $-\Delta_{S^{n-1}} + \left(\frac{n-2}2\right)^2$. Its spectrum 
is 
\[ \{(\frac{n-2}2 + k)^2:=\sigma_k^2 : k=0,1,\dots \} \] 
and the corresponding eigenspace is spanned by harmonic polynomials. 
The equation becomes
\[ u^k_{tt} - \sigma_k^2 u^k = e^{-\frac{n+2}{2}t} g^k  \]
and the estimate \eqref{carleman1}  follows from (including an additional linear term into $h$ without changing the notation) 
\[  \int e^{2h} (1+h'') (|u^k_t|^2+ (1+k^2 ) |u^k|^2 +|h'|^2 |u^k|^2) dt 
\le C \int e^{2h-4t} |g^k|^2  dt \]     
Since $\partial_t^2 -\sigma_k^2 = (\partial_t -\sigma_k)(\partial_t+\sigma_k) $, 
the claim follows once we prove the elliptic estimate 
\begin{equation} \label{factor1}   \int e^{2h}( |u'|^2 + (\tau^2+\sigma^2) |u|^2 ) dt \le  C \int e^{2h}  |(\partial_t -\sigma) u |^2 dt \end{equation} 
and the commutator type estimate 
\begin{equation}   \int e^{2h} (1+h'')|u|^2  dt \le  c \int e^{2h}  |(\partial_t +\sigma) u |^2 dt \end{equation}

In the first case we multiply 
\[  u_t -\sigma u = g \] 
by $e^{2h}u$ and integrate. Then 
\[ 
\frac12 (\tau+\sigma) \Vert e^{h} u \Vert^2 
\le  \int e^{2h} (h'+ \sigma) u^2 dt = - \int e^{2h} u  g dt 
\le  \Vert e^{h} g\Vert\, \Vert e^{h} u \Vert
\] 
together with using the equation to bound $u_t$  implies \eqref{factor1}. For the second estimate we define $v = e^{h} u$, multiply 
\[  v' - (h'-\sigma) v = e^hg \] 
by $(h'-\sigma) v$ and obtain 
\[ 
\Vert (h'-\sigma) v \Vert^2 + \frac12 \int h'' v^2 dt = -\int e^h g (h'-\sigma)v dt. 
\]  
The estimate follows by an application of the Cauchy-Schwarz inequality.

\eproof

Besides the Carleman estimate, we also need an interior estimate (Caccioppoli-type estimate) for the Lam\'e system \eqref{1.2}. For fixed $a_3<a_1<a_2<a_4$, there exists a constant $C_1$ such that
\begin{equation}\label{3.1}
\int_{a_1r<|x|<a_2r}||x|^{|\alpha|}D^{\alpha} u|^2+||x|^{|\alpha|+1}D^{\alpha} p|^2dx\le
C_1\int_{a_3r<|x|<a_4r}|u|^2dx,\quad|\alpha|\le 1
\end{equation}
for all sufficiently small $r$. Estimate \eqref{3.1} can be found in Lemma 3.1 of \cite{lin6}.

We are now ready to prove Theorem~\ref{thm1.1} and \ref{thm1.2}. Let us define the cut-off function  $\chi(x) \in
C^{\infty}_0 ({\mathbb R}^n\backslash {\{0\}})$ such that
$$
\chi(x)=\begin{cases}
0\quad\text{if}\quad
|x|\leq r/3,\\
1\quad\text{in}\quad 5r/12\leq|x|\leq 5\tilde R/4,\\
0\quad\text{if}\quad 3\tilde R/2\leq |x| ,
\end{cases}
$$
where $\tilde R$ is a small number that will be chosen later and $r\ll \tilde R$. Denote $\tilde u=\chi u$ and $\tilde p=\chi p$. Then it follows from \eqref{up} that $\tilde u$ and $\tilde p$ satisfy
\begin{equation}\label{eq3-2}
\Delta\tilde u+\nabla(a\tilde p)=(\nabla^2\chi)u+2\nabla\chi\cdot\nabla u+(\nabla\chi) ap-\chi G:=F
\end{equation}
and
\begin{equation}\label{eq3-3}
\Delta\tilde p+\mbox{div}(\chi G)=(\nabla^2\chi)p+2\nabla\chi\cdot\nabla p+(\nabla\chi) G:=H.
\end{equation}
Applying \eqref{carleman1} to \eqref{eq3-2} with $u=\tilde u$, $f=a\tilde p$, $g=F$ yields
\begin{equation}\label{eq3-5}
\begin{aligned}
&\tau \Vert   |x|^{-2} (1+h'')^\frac12 e^h \tilde u \Vert + \Vert  |x|^{-1} (1+h'')^\frac12 e^h  \nabla\tilde u \Vert\\
&\lesssim \tau  \Vert |x|^{-1}  e^h a\tilde p \Vert +\Vert e^h F\Vert\\
&\le C (\tau  \Vert |x|^{-1}  e^h \tilde p \Vert +\Vert e^h F\Vert),
\end{aligned}
\end{equation}
where $C=C(n,M_0,\delta_0)$. Replacing $h$ by $h-t-\frac 12\ln(1+h'')$ in \eqref{carleman1} and applying the new estimate to \eqref{eq3-3}, we have that
\begin{equation}\label{eq3-6}
\tau \Vert   |x|^{-1} e^h \tilde p \Vert + \Vert  e^h  \nabla\tilde p \Vert\lesssim \tau  \Vert (1+h'')^{-\frac 12} e^h\chi G \Vert +\Vert |x|(1+h'')^{-\frac 12}e^h H\Vert.
\end{equation}
Now, $K\times$\eqref{eq3-6}$+$\eqref{eq3-5} gives
\begin{equation}\label{eq3-8}
\begin{aligned}
&\tau \Vert   |x|^{-2} (1+h'')^\frac12 e^h \tilde u \Vert + \Vert  |x|^{-1} (1+h'')^\frac12 e^h  \nabla\tilde u \Vert+K\tau \Vert   |x|^{-1} e^h \tilde p \Vert\\
&\le C(\tau  \Vert |x|^{-1}  e^h \tilde p \Vert +\Vert e^h F\Vert+ K\tau  \Vert (1+h'')^{-\frac 12} e^h\chi G \Vert +K\Vert |x|(1+h'')^{-\frac 12}e^h H\Vert).
\end{aligned}
\end{equation}
We then choose $K\ge C$ and $\tilde R=R(n,M_0,\delta_0)$ satisfying
\[
CK\tau(1+h'')^{-\frac 12}\le |x|^{-1}(1+h'')^{\frac 12}
\]
for all $|x|\le 3R/2$ since \eqref{eq333} holds. Consequently, we obtain from \eqref{eq3-8} that
\begin{equation}\label{310}
\begin{aligned}
&\tau \Vert   |x|^{-2} (1+h'')^\frac12 e^h u \Vert_{\{5r/12\le|x|\le 5R/4\}}\le C\Vert e^h F\Vert_{\{r/3\le|x|\le 5r/12\}\cup\{5R/4\le|x|\le 3R/2\}}\\
&+C \tau  \Vert (1+h'')^{-\frac 12} e^hG \Vert_{\{r/3\le|x|\le 5r/12\}\cup\{5R/4\le|x|\le 3R/2\}}\\
&+C\Vert |x|(1+h'')^{-\frac 12}e^h H\Vert_{\{r/3\le|x|\le 5r/12\}\cup\{5R/4\le|x|\le 3R/2\}}:={RHS}.
\end{aligned}
\end{equation}
Here and after, we use $\|\cdot\|_A$ to denote the $L^2$ norm over the region $A$.

In view of \eqref{3.1}, we can easily derive that
\begin{equation}\label{320}
RHS\le C\tau e^{\tilde h(r/3)}(r/3)^{-2}\Vert u\Vert_{\{r/4\le|x|\le r/2\}}+C\tau e^{\tilde h(5R/4)}(5R/4)^{-2}\Vert u\Vert_{\{R/2\le|x|\le 2R\}},
\end{equation}
where we denote $\tilde h(a)=h(-\ln a)$. Now we choose $\tau=\tau_0$ such that
\begin{equation}\label{321}
Ce^{\tilde h(5R/4)}(5R/4)^{-2}\Vert u\Vert_{\{R/2\le|x|\le 2R\}}\le\frac 12 R^{-2}e^{\tilde h(R)}\Vert u\Vert_{\{2R/3\le|x|\le R\}}.
\end{equation}
More precisely we  choose from now on
 \[ \tau_0\sim
\ln\Big( \frac{\Vert u\Vert_{\{2R/3\le|x|\le R\}}}{\Vert u\Vert_{\{R/2\le|x|\le 2R\}}}\Big).
\]
so that \eqref{321} is satisfied.
Combining \eqref{310}, \eqref{320}, \eqref{321} yields
\begin{equation}\label{322}
\Vert |x|^{-2} e^h u \Vert_{\{5r/12\le|x|\le 5R/4\}}\le Ce^{\tilde h(r/3)}(r/3)^{-2}\Vert u\Vert_{\{r/4\le|x|\le r/2\}}.
\end{equation}
The estimate implies that
\[
\Vert u\Vert_{\{|x|\le r\}}\ge Ce^{\tilde h(R)}\Vert u \Vert_{\{2R/3\le|x|\le R\}}(r/R)^2e^{-\tilde h(r/3)}\ge C r^{m},
\]
which establishes Theorem~\ref{thm1.1}. Next, adding $e^{h(r/2)}(r/2)^{-2}\Vert u\Vert_{\{|x|\le r/2\}}$ to both sides of \eqref{322} gives
\[
\begin{aligned}
e^{\tilde h(r)}r^{-2}\|u\|_{\{|x|\le r\}}&\le e^{\tilde h(r)}r^{-2}\|u\|_{\{|x|\le r/2\}}+e^{\tilde h(r)}r^{-2}\|u\|_{\{r/2\le|x|\le r\}}\\
&\le Ce^{\tilde h(r/3)}(r/3)^{-2}\|u\|_{\{|x|\le r/2\}},
\end{aligned}
\]
which leads to Theorem~\ref{thm1.2}.

\section{Global doubling inequalities}\label{sec4}
\setcounter{equation}{0}

In the previous section, we have proved local doubling inequalities. Nonetheless, global doubling inequalities are more suitable for inverse problems (for example, see \cite{amr}). In this section we derive global doubling inequalities along the
lines in \cite{amrv}. For brevity, we will not give detailed
arguments here. We refer to \cite{amrv} for detailed proofs. To
begin, we give the definition of Lipschitz boundary.
\begin{definition}\label{deflip}
We say that the boundary $\partial\Omega$ is of Lipschitz class with
constants $r_0$ and $L_0$, if, for any $x_0\in\partial\Omega$, there
exists a rigid transformation of coordinates under which $x_0=0$ and
$$\Omega\cap B_{r_0}(0)=\{x\in B_{r_0}(0): x_n>\psi(x')\},$$ where
$x=(x',x_n)$ with $x'\in\R^{n-1}, x_n\in\R$ and $\psi$ is a
Lipschitz continuous function on $B_{r_0}(0)\subset\R^{n-1}$
satisfying $\psi(0)=0$ and
$$\|\psi\|_{C^{0,1}(B_{r_0}(0))}\le L_0r_0.$$
\end{definition}

Let us denote
$\Omega_d=\{x\in\Omega:\text{dist}(x,\partial\Omega)>d\}$. Using
three-ball inequalities proved in \cite{lin5} or \cite{lin6}, one can prove the following theorem (see
\cite{amr}, \cite{amrv}).
\begin{theorem}{\rm\cite[Theorem~3.2]{amrv}}\label{propaga}
Let $\partial\Omega$ be of Lipschitz class with constants $r_0$,
$L_0$, and $\lambda$, $\mu$ satisfy \eqref{1.1} , $u\in H^1_{loc}({\Omega})$ be a nontrivial solution to \eqref{1.2}.
Then for every $\sigma>0$ and for every
$x\in\Omega_{\frac{4\sigma}{\theta}}$, we have
$$\int_{B_{\sigma}(x)}|u|^2dx\ge C_{\sigma}\int_{\Omega}|u|^2dx,
$$
where $0<{\theta}<1$
depends on $n$, $\delta_0$, $M_0$ only and $C_{\sigma}$
depends on $n$, $\delta_0$, $M_0$, $r_0$, $L_0$, $|\Omega|$,
$\|u\|_{H^{1/2}(\Omega)}/\|u\|_{L^2(\Omega)}$, and $\sigma$.
\end{theorem}

We now ready to state global doubling inequalities. To describe the theorem, we introduce more notations. Instead of the
strong ellipticity, we say that Lam\'e coefficients $\lambda$, $\mu$
satisfy the strong convexity condition if
\begin{equation}\label{convex}
\mu(x)\ge\tilde\delta_0>0,\quad
2\mu(x)+n\lambda(x)\ge\tilde\delta_0\quad\forall\quad x\in \Omega.
\end{equation}
It is known that the strong convexity implies the strong
ellipticity. Let $\varphi\in L^2(\partial\Omega,\R^n)$ be a vector
field satisfying the compatibility condition
$$\int_{\partial\Omega}\varphi\cdot rds=0$$ for every infinitesimal
rigid displacement $r$, that is, $r=c+Wx$, where $c$ is a constant
vector and $W$ is a skew $n\times n$ matrix. Consider the boundary
value problem:
\begin{equation}\label{bvp}
\begin{cases}
\div(\mu(\nabla u+(\nabla u)^T))+\nabla(\lambda\div
u)=0\quad\text{in}\quad\Omega,\\
(\mu(\nabla u+(\nabla u)^T)+(\lambda\div
uI_n))\nu=\varphi\quad\text{on}\quad\partial\Omega,
\end{cases}
\end{equation}
where $I_n$ is the $n\times n$ identity matrix, $\nu$ is the unit
outer normal to $\partial\Omega$, and $\varphi$ satisfies the
compatibility condition. In order to ensure the uniqueness of the
solution to \eqref{bvp}, we assume the following normalization
conditions:
\begin{equation}\label{normal}
\int_{\Omega}u dx=0,\quad\int_{\Omega}(\nabla u-(\nabla u)^T)dx=0.
\end{equation}
\begin{theorem}{\rm\cite[Theorem~3.7]{amrv}}\label{gdouble}
Let $\partial\Omega$ be of Lipschitz class with constants $r_0$,
$L_0$, and $\lambda$, $\mu$ satisfy \eqref{convex} , the second condition of \eqref{1.1}. If
$u\in H^1(\Omega,\R^n)$ is the weak solution to \eqref{bvp}
satisfying the normalization condition \eqref{normal}. Then there
exists a constant $0<\vartheta<1$, only depending on $n$,
$\tilde\delta_0$, $M_0$, such that for every $\bar{r}>0$ and for
every $x_0\in\Omega_{\bar r}$, we have
$$
\int_{B_{2r}(x_0)}|u|^2dx\le C\int_{B_r(x_0)}|u|^2dx
$$ for every $r$ with $0<r\le\frac{\vartheta}{2}\bar r$, where $C$
depends on $n$, $\tilde\delta_0$, $r_0$, $L_0$, $|\Omega|$, $\bar
r$, and
$\|\varphi\|_{H^{-1/2}(\partial\Omega)}/\|\varphi\|_{H^{-1}(\partial\Omega)}$.
\end{theorem}

\renewcommand{\theequation}{A.\arabic{equation}}
\renewcommand{\thelemma}{A.\arabic{lemma}}

\setcounter{equation}{0}  
\section*{Appendix}  

In this appendix, we would like to construct a weight function $h$ satisfying the conditions described in Section~\ref{sec3}. 
Let $\tau\in \mathbf{N}+\frac{5}{4}\gg 1$ and define $a=2\ln \tau$. We choose
$$
h''(t)=\delta \tau e^{-t/2},
$$
where $\delta>0$ is sufficiently small. We then set
$$
h'(t)=\tau-2\delta\tau e^{-t/2}
$$
and
$$
h(t)=\tau t+4\delta\tau e^{-t/2}.
$$

It is clear that $h$ is convex and $h'$ satisfies the first condition of \eqref{weight}. To verify the second condition of \eqref{weight}, we observe that
$\tau e^{-t/2}\leq 1$ if $t\geq 2\ln \tau(=a)$ and $\tau e^{-t/2}\geq 1$ if $t\leq 2\ln \tau$. So, for $t\leq a$, we have $h''(t)\ge\delta \tau e^{-t/2}\ge C_\delta(1+\tau e^{-t/2})$ for some $C_\delta>0$. Next, for $a< t$, we can see that $\tau-2\delta\leq h'(t)\leq \tau$, then $\mbox{dist}(2h',{\mathbb Z})\geq \frac{1}{2}-4\delta\ge C(1+\tau e^{-t/2})$ holds for some absolute constant $C>0$ provided $\delta\le\frac{1}{16}$. 

To check \eqref{eq333}, as we noted above, if $t\le 2\ln\tau$, then
\[
1+h''(-\ln|x|)\ge 1+\delta\tau \sqrt{|x|}\ge \delta\tau |x|
\]
for $|x|<1$. On the other hand, for $t\ge 2\ln\tau$, we have
\[
1+h''(-\ln|x|)\ge 1\ge\tau e^{\ln|x|/2}=\tau\sqrt{|x|}\ge \tau |x|.
\]
Finally, let us define $\tilde{h}=h-t-\frac 12\ln(1+h'')$, then we have
$$
\tilde{h}''(t)=\delta \tau e^{-t/2}-\frac{1}{8}\delta \tau e^{-t/2}(1+\delta \tau e^{-t/2})^{-1}+\frac{1}{8}\delta^2 \tau^2 e^{-t}(1+\delta \tau e^{-t/2})^{-2}.
$$
We choose
$$
\tilde{h}'(t)=\tau-2\delta\tau e^{-t/2}-1+\frac{1}{4}\delta \tau e^{-t/2}(1+\delta \tau e^{-t/2})^{-1}
$$
and
$$
\tilde{h}(t)=\tau t+4\delta\tau e^{-t/2}-t-\frac 12\ln(1+\delta \tau e^{-t/2}).
$$
The same arguments imply that $\tilde{h}$ satisfies the required conditions provided $\delta$ is small.

\end{document}